\documentclass{article}    

\usepackage{graphicx}
\title{\textbf{Non-trivial 1-classes in the homology of the real moduli spaces $\overline{M}_{0,n}$ and related structures}}  
\author{Gefry Barad\\{\small Department of Mathematics, Johns Hopkins University}\\{\footnotesize Baltimore MD 21218, gefry@math.jhu.edu, http://www.math.jhu.edu/\~{}gefry}}         
\begin{document}
\label{}          
\maketitle                 
\begin{abstract}{{\normalsize We give lower bounds for the rank of the first homology group of the real points of the Deligne-Mumford-Knudsen compactification of stable n-pointed curves of genus 0,which coincides with the Chow quotient $(RP^1)^n//PGL(2,R)$.~\cite{ka}. The study is a diagrammatic version of Loday's Theorems ~\cite{loday}.Also, it has connections with spectral geometry and combinatorics.}}\\\textbf{A.M.S. Classification:}14P25 (14F35, 14J10, 57M15, 20F99) 
                                                                                                                                   \end{abstract}

{\Large \section{Discrete Spectral warm-up} \label{bella, morava}
We built a sequence of 2 dimensional complexes. Its study ,
inspired by spectral graph theory, can guide the study of the real
moduli spaces $\overline{M}_{0,n}$. Also, these complexes give a
counterexample to a tempting conjecture.

\vspace{6pt}
\par\noindent

1. Let G(n) be the Schreier coset graph of S(n), the symmetric
group, modulo D(n), the dihedral group, with respect to the
following set T of generators. T is a set of $n(n-3)/2$
involutions from $S_n$. For every pair (k, l), $0<k<l<n+1$, where
(k,l) is different from (1, n), we define d:=d(k, l)(x)=l+k-x , if
x is between k and l. Otherwise, x is a fixed point. \vspace{10pt}

The vertices of G(n) are the right coset classes S(n)/D(n). There
is an edge between D(n)g and D(n)h , if there is t in T such that
D(n)g=D(n)ht.

\vspace{10pt}

Any regular graph ( a graph where any vertex has the same degree)
is a Schreier coset graph of a group G, modulo a subgroup H, with
respect to a set of generators T.~\cite{schreier} \vspace{10pt}

Let B(n) the following 2 dimensional cellular complex. Its
1-skeleton is given by G(n). For any u and v in T, such that u and
v commutes, there is a face around D(n)a, D(n)av, D(n)avu, and
D(n)au.

As we will prove, rank $H_1$(B(n)) = rank
$H_1(\overline{M_{0,n}}(R))$. So, the study of B(n) is important.
The first barycentric subdivision transforms any cellular complex
into a simplicial one.

\vspace{10pt}

2.\textit{Definition 1.} For any 2 dimensional simplicial complex,
the link of a vertex v is the graph L(v) whose vertices are the
edges incident to v, and whose edges are the faces between 2
incident edges to v. So, the graph is the intersection between the
cone from v and a plane.

\vspace{10pt}

The Laplacian of a simple finite connected graph is the matrix of
entries a(i,j)/$\sqrt{(d(i)d(j)}$, where the rows and columns are
indexed over vertices of the graph. d(v) is the degree of a vertex
v, i.e the number of vertices joined with v. a(i,j) is 1 if and
only if there is an edge between i and j.~\cite{chung}.

\vspace{10pt}

 Why do we need this graph? The cohomology with real
coefficients of a compact simplicial complex is zero IF:
\par\noindent
a. [Ballmann]: for every vertex v, the first positive eigenvalue
of the Laplacian on L(v) is greater than 1/2. See also ~\cite{s}.
\par\noindent
or IF:
\par\noindent
b. [Zuk]: for any 2 vertices a and b of the complex, joined by an
edge, the average of their first positive eigenvalue of the
Laplacians on L(a) and L(b) is greater than 1/2.~\cite{zuk}.

\vspace{10pt}

Let us formulate \textbf{the tempting conjecture: } The first
cohomology with real coefficients is zero, if the average of first
positive eigenvalue of  any 3 vertices of the same triangle is
greater than 1/2. This conjecture is not true. A counterexample is
given by the spaces B(n). As we will see, their first cohomology
is always different from zero, even if the conditions from the
statement of this conjecture are satisfied for n greater than 70!
\vspace{10pt}

The first barycentric subdivision of B(n) has 3 types of vertices.
The first one is a center of a face. Its link is an 8-gon. The
first positive eigenvalue is 1 - sin($\pi$/4). The second one is a
center of an edge. Its link is a bipartite graph and its first
positive eigenvalue is 1. A wonderful treatment of this
theoretical and computational subject is give in ~\cite{graph}.
Every vertex from G(n) has the same isomorphic link, because of
the regularity of the Schreier coset graphs. The last type of
vertex is a vertex from the old G(n). The first positive
eigenvalue of its link, $\lambda$(n) was computed by Matlab 6.1
using a C++ program, available at the author's web page.
$\lambda$(n) is an increasing function of n.

\vspace{10pt}

$\lambda$(6)=0.1888     $\lambda$(7)=0.1891 $\lambda$(50)=0.2070
$\lambda$(70)=0.2084

\vspace{10pt}

There is the following explanation of the increasing behavior of
$\lambda$(L(n)),(L(n)= the link of a vertex from G(n)).
\vspace{10pt}

 L(n) is the graph whose vertices are the diagonals of an n-gon. 2
 vertices are joined by an edge if the interiors of the diagonals
 do not intersect. L(n) has a D(n) symmetry. Also, it is included
 in L(n+1).If K(n) is the complete graph of n vertices, there is a covering map from the cartesian product L(n) X K(n+1) to
 L(n+1) . So, any eigenvalue of the latter graph is an eigenvalue
 of the former graph , too (see Chung). Also, the eigenvalues of
 a cartesian product are sums of separate eigenvalues, because the Laplacian of a cartesian product is a tensor product of separate matrices. Zero is
 always an eigenvalue, so the smallest eigenvalue of the cartesian product is the minimum between the smallest eigenvalue of L(n) and K(n). In this case, it is the smallest positive eigenvalue of L(n). We've got the desired inequality between
 2 consecutive $\lambda$(n)'s.

Note: applying Zuk's Theorem to hypergraphs.~\cite{zuk}, we get
that any closed path in G(n) is a sum of closed paths of length 4.
We did not find an equivalent statement in algebraic geometry.

\section{The first homology group of $\overline{M_{0,n}}(R)$}

In their Quantum Invariants and Knot Theory books and articles,
Turaev, Kassel and Dror bar-Nathan use "graphical proofs" and
diagrammatic calculus, for algebraic or categorical data. There is
an intricate relation between formal algebraic structures and
concrete geometrical objects. Even if it is elementary, it is
useful to describe the path towards $H_1( \overline{M_{0,n}}(R))$
because this is a way to understand the proofs and the formalism
of Loday and the appearance of these spaces as Operads.

\vspace{10pt}

To keep track of the first homology group,  we will use the
cellular decomposition of M=$\overline{M}_{0,n}(R)$, after
Devadoss . We apply the dual block complex, ~\cite{mu}pp.380, and
restrict to level two of its filtration.

\vspace{10pt}

Let $K_{n-1}$ be the n-3-convex polytope whose partial order set
of its faces is isomorphic with the partial order set of an n-gon
with a couple of non-intersecting diagonals.  The partial orders
are given by inclusions. $K_{n-1}$ is called the associahedron.
There is such a convex polytope [Ziegler-Lectures on Polytope
p.310]!  $K_{n-1}$ has the acyclic carrier property , so we can
apply the classical theorems from [Munkres pp 225]. The
codimension k faces of the associahedron are indexed by n-gons
with k non-intersecting diagonals.

\vspace{10pt}

DIAG is the following set of $n(n-3)/2$ involutions from $S_n$.
For every pair (k, l), $0<k<l<n+1$, where (k,l) is different from
(1, n), we define d:=d(k, l)(x)=l+k-x , if x is between k and l.
Otherwise, x is a fixed point. Let P be a fixed n-gon, with edges
labeled 1,2,3...n.  For every diagonal of P we can associate an
element of DIAG in the following way:  any diagonal determines a
partition of 1,2...n.  Take the one which doesn't contain n:  it's
between 2 numbers, k and l.  Then the associated d will be d(k,l),
and we say that d(k,l) is supported by the diagonal of the n-gon
P.  \textbf{Throughout the paper, the word "diagonal" means a
diagonal of the n-gon, or the involution carried by the diagonal.}

\vspace{10pt}

Take $n!$ copies of $K_{n-1}$.  For every permutation of $S_n$, label the edges of the n-gon with $\sigma{(1)}, \sigma{(2)}..... \sigma{(k)}, ....\sigma{(n)}$.  So the codimension k faces are labeled by decorated n-gons with k non-intersecting diagonals.  Now we build our space $\overline{M}_{0,n}(R)$.  Two codimension k faces of different $K_{n-1}$'s are identified (glued) if the permutations $\sigma_1$ and $\sigma_2$ which color the edges of the n-gons satisfy the following condition ``flip" or gluing condition:  there are $d_1, d_2, ....,d_i$ couple of elements of DIAG, supported by the diagonals of the second face, such that $\sigma_{1} = \sigma_{2}\circ d_1\circ d_2\circ...\circ d_i$. ( $\circ$ means composition of functions ).
\vspace{12pt}

The top dimensional faces (without diagonals) are identified by the action of $D_n$, the dihedral group.  So we can begin with $(n-1)!$/2 copies of $K_{n-1}$, indexed over $S_n$/$D_n$.  Two codimension k faces are identified if their classes modulo the dihedral group $D_n$ contain 2 permutations which satisfy the flip condition from the previous paragraph.

\vspace{10pt}

\textbf{The Homology} is encoded in the gluing process above.
M(our moduli space) is a smooth compact (n-3)-manifold,
non-orientable if the dimension is higher than 1.

Recall the construction of the dual block complex.  Let X be a compact homology n-manifold.  Let sdX be the first barycentric subdivision.  The simplices of sdX are $[\bar{a}, \bar{b}...\bar{z}]$ where ${a}\supset{b}.....\supset{z}$ and $\bar{\sigma}$ is the barycenter of $\sigma$ , a simplex of X.

\vspace{10pt}

Given a simplex $\sigma$, D($\sigma$), the block of $\sigma$, is the union of the open simplices of sdX, where $\sigma$ is the final vertex; i.e. "$\sigma$ = z" in the notation above.  We have dim $\sigma$ + dim D($\sigma$) = dim manifold.

\vspace{10pt}

2.0.Let $X_p$ be the dual p-skeleton of X, the union of all
D($\sigma$) such that dim D($\sigma$) is smaller or equal to p.
$\textbf{D}_p$ = $H_p$($X_p$ , $X_{p-1}$).  The boundary operator
is the boundary operator in the exact sequence of the triple
($X_p$ , $X_{p-1}$ , $X_{p-2}$). $\textbf{D}_p$ = $H_p$($X_p$ ,
$X_{p-1}$) is the free abelian group generated by the blocks of
dimension p of X, arbitrarily oriented, i.e. the blocks from the
previous statement label the elements of a basis. How can we deal
with the boundary operator of the dual block complex?  Fortunately
it has a nice geometric meaning.  The dual blocks form a CW
decomposition of X.  Let a and b be 2 cells, dim(a)=p and
dim(b)=p-1. $\partial{a}$ is a formal sum of p-1 cells, with
integer coefficients.  The coefficient of b is given by "the
incidence coefficient", which is the degree of a map that sends
the boundary of a (i.e. $S^{p-1}$) to a bouquet of p-1 spheres =
$X_{p-1}$/$X_{p-2}$, and then projected to a p-1 sphere.  So the
boundary map shows how the boundary of the cell is patched by p-1
spheres.

\vspace{10pt}

2.1 We would like to apply the previous settings to X =
$\overline{M}_{0,n}(R)$. The associahedra give a cellular, not a
simplicial decomposition of X. We have to take the first
barycentric subdivision of X. Fortunately, the barycenters of the
faces of the associahedra are already labelled by n-gons with a
couple of diagonals, where the edges of the n-gons bear a
permutation.

\vspace{10pt}

A maximal simplex in X is a sequence of n-3 barycenters, i.e a
chain of length n-3, which joins 2 vertices at the distance n-2 ,
in the graph G(n) from the previous section. Its dual block, of
dimension 0, is its barycenter, which can be labelled by the
simplex itself. In the notation from 2.0 section, $\textbf{D}_{0}$
is the free abelian group generated by these simplices.

\vspace{7pt}

A codimension 1 simplex of a simplex above is a face of the
simplex above. Its dual block is a segment between 2 maximal
simplices  which share the same face (X is a manifold !).
$\textbf{D}_{1}$ is the free abelian group generated by these
segments, arbitrarily oriented.

\vspace{7pt}

There are 2 types of segments: the segments inside the same
associahedron, and the segments between 2 different associahedra.

Similarly, $\textbf{D}_{2}$ is the free abelian group generated by
4 segments which form the boundary of a 2-cell, arbitrarily
oriented. Any dim 2 block is shared by exactly 4 dim 1 blocks and
exactly 4 dim 0 blocks, thereby building a structure similar with
the structure of 4 cubes in 3 dimensions.

\vspace{6pt}

The boundary morphisms between D's are "normal": a segment goes to
the difference of its vertices, and a 2-face goes to a sum of
edges, correlated by signs.

\vspace{8pt}

$H_1(\overline{M}_{0,n}(R))$ =$Ker(\partial_{1})/Im(\partial_{2})$

\vspace{7pt} The barycenters of the maximal simplices above and
the edges among them form a new graph, called GG(n). There is a
following pictorial transformation between G(n) and GG(n):
\vspace{5pt} -the vertices of G(n) become circles. Between these
circles, instead of one edge e , there are m(e) edges. m(e) is the
following number: any edge e is decorated by a diagonal of the
n-gon. m(e) is the number of diagonals which do not intersect e.
The tips and the tails of these edges are inside the circles and
there are connected by a system of pipes given by the barycentric
subdivision of the associahedra.

\includegraphics [scale=2.0]{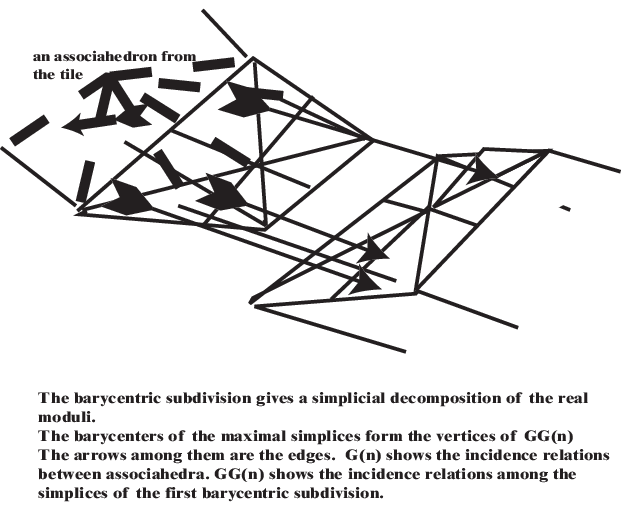}

An element of $D_{1}$ is a formal sum of edges, decorated with
reals. We can decorate the edges of with any algebraic structure
V. Between 2 vertices of G(n) there is a vector from the tensor
product of m(e) copies of V.

Theorem 1. There is an  isomorphism between the first homology
groups \textit{with real coefficients} of B(n) and
$\overline{M}_{0,n}$. It is given by the following function F: If
the m(e) edges (between circles) are decorated by a couple of real
numbers in GG(n), we associate their sum to the edge e in G(n).
Using the fact that the associahedron is homeomorphic with a
closed ball, it is easy to prove the statement above. It is true
only for real coefficients !There are 3 questions to be answered
in this proof: \vspace{10pt}

1. If the definition is independent of the numbers (it depends
only on the homology class). the image of a boundary is a
boundary. By a boundary, we mean a formal sum of 4 edges, colored
by + or -1, according to their orientation , which are geometric
boundaries of 2-blocks etc.

2. If we can define an inverse of F. We can define a function G
from the homology of $B(n)$ to $H_{1}(M)$, defined on generators
in the following way. : the edge e colored by number x goes to the
sum of the m(e) edges from GG(n),colored by the same number x/m(e)
- a kind of trace-diagonal process. We do not have any
obstructions: the associahedron is aspherical, so it is possible
to assign numbers to the edges inside them , such that the result
is a flux in GG(n). ( a flux is an assignment of numbers to edges
such that the sum for every vertex is zero, see section 4). A
physics of this process is given by the concept of "pressure": the
pressure of a gas inside the associahedron is zero. it is just a
distribution of pressures given by the numbers from the external
faces.

3. If they are inverse to each other. It is easy, to see that this
is the case.

\vspace{12pt}

\section{Koszulness and Koszul Duality for Operads.\\\
Combinatorics. The beginning of a project.}

If, instead of real numbers, we decorate the graphs above with
vector spaces, and "arrows" means morphisms (we already fixed an
arbitrarily orientation for edges), our result is exactly the
Theorem 4.3 of Loday ~\cite{graph}: the simplex and the Stasheff
quadratic operads are Koszul.

\vspace{5pt} The 2 operads above are Koszul dual to each other. A
natural question is: what is the dual diagrammatic calculus?

\vspace{5pt}

An associahedron can be divided in standard simplices using the
barycentric subdivision. Is it possible to divide a standard
simplex in a couple of associahedra? For example, a 2-simplex can
be divided in 9 pentagons.

\includegraphics [scale=2.0]{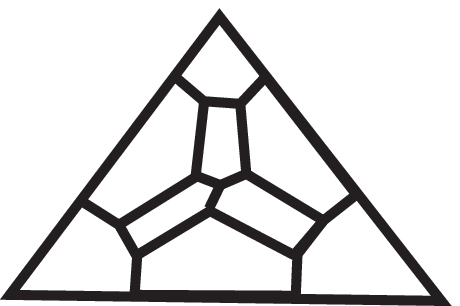}

Maybe an inductive study can show that it is true for any n. Is it
the right structure for the dual picture?

In a very recent paper by Stanley ~\cite{stanley} the "right"
picture showed up: \textit{A simplex can be divided in chambers,
such that these chambers are products of lower dimensional
simplices, and the partial order set , under inclusion, of the
faces of the chambers are isomorphic with the associahedron.}We
call it "the associahedron internal to an n-simplex". Instead of
"associahedra" we can use standard simplices, glued together as
the corresponding associahedra are glued together. This is just
the beginning of the picture. An axiomatic treatment is needed to
complete the duality. \vspace{3pt} Also, it is worth to study the
universality or the freeness of the above structures, as the role
played by tangles for braided categories.

The subject above opened the following questions:
\par\noindent
\vspace{10pt}
1. Is it possible to divide an n-simplex in a couple
of associahedra. Is it possible to apply induction?
\par\noindent
2. Does the Stanley construction give rise to a simplicial complex
equivalent with the real moduli spaces? A study of the
automorphism group of Stanley's construction can show how can we
glue the faces of the simplices. This simplicial decomposition can
be useful in homological computations.
\par\noindent
3. Is it possible to build an operad based on simplices, following
the ideas above ? Lower dimensional faces of associahedra are
products of lower dimensional associahedra. The chambers of the
decomposition of an n-simplex are products of lower simplices. The
natural boundary or co-boundary maps among chambers can give this
structure.
\par\noindent

\section{Non-trivial 1-classes in the homology of the real moduli spaces $\overline{M}_{0,n}$}
\vspace{2pt}
Because of the homology isomorphism above, it is
enough to study $B(n)$, to get the desired results of this
section. A partial study of these spaces is given in
~\cite{barad}.

The notes below is a study of real vector spaces or subspaces
generated by formal sums of oriented edges in G(n).
 We found 4 special types of cycles in B(n):

1. Cycles given by k non-intersecting diagonals in the n-gon, and
by k numbers of zero sum. (picture 1 for k=3).

\includegraphics [scale=2.0]{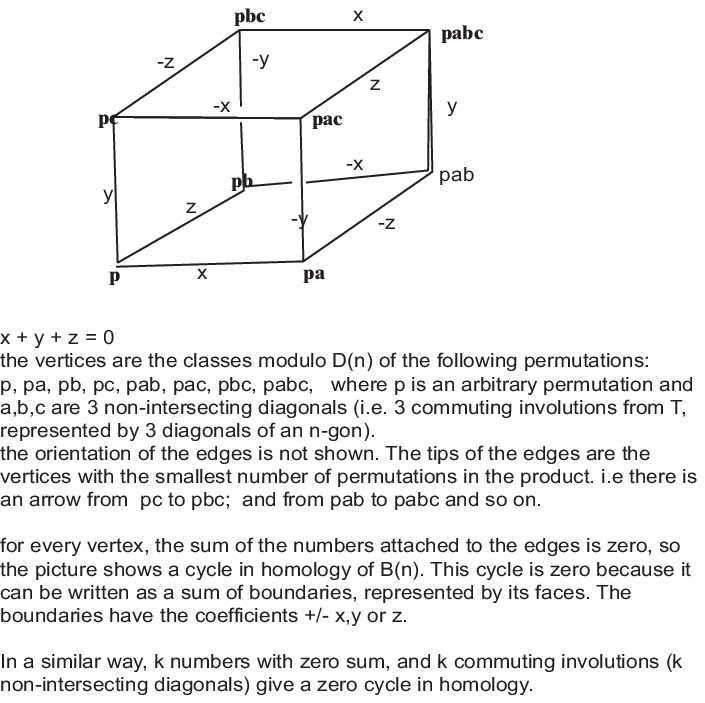}

These cycles are boundaries, they are zero in homology. A cycle of
this type is the 1-skeleton of a k-dimensional cube. For every k
numbers $a_{i}$ with zero sum, we can decorate the vertices of the
cycle with + or - $a_{i}$ in such a way we get a cycle. We can
write the 1-skeleton like a sum of cycles of length 4, so it is
zero in homology.

2. Cycles given by regular k-gons. (picture 2 for k=4)

\includegraphics [scale=2.0]{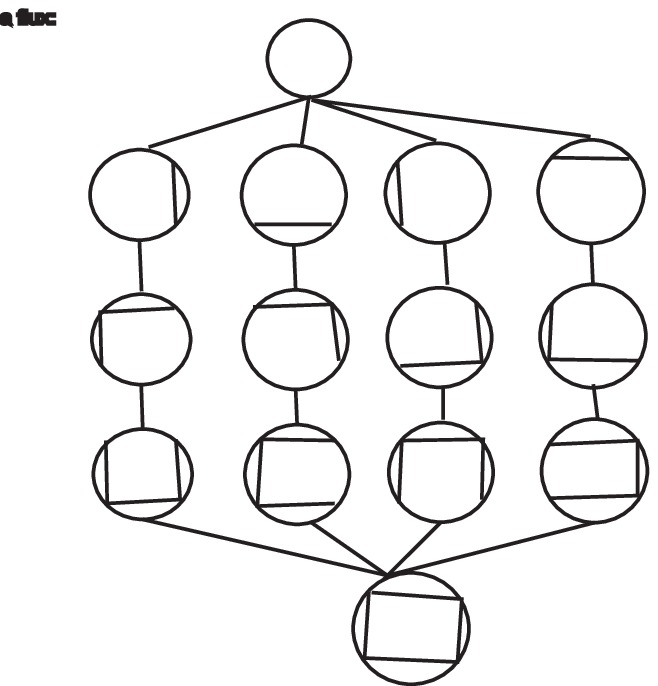}

3. Cycles given by Reiner's duality ~\cite{vic}: Let a and b 2
vertices of G(n) at the distance n-2. The vertices which are on
the paths of length n-2 form a self-dual lattice. So, it is enough
to color the edges of G(n) with numbers, only for the first half
of the lattice. The remaining part is the anti-symmetric.

    We did not find a way to decide if the previous
    2 types of cycles are zero or not in homology.

 4. Cycles give by the graph $G(n,k)$, the subgraph of $G(n)$ induced by the vertices of distance at most k from the identity vertex.

 These are cycles of "finite type".
 For k=2, the computation of the dimension M of the space of cycles
 was made in ~\cite{barad}. These cycles survive in the whole
 homology of B(n), they cannot be cancelled by external
 boundaries.

 Let us mention the formulas for M, and \textbf{{\large a sketch  of the proof.}
}
\begin{displaymath}
\texttt{} \hspace{0.5em}M = \left\{
\begin{array}{ll}
1 & \textrm{if $n=4$}\\
4 & \textrm{if $n=5$}\\
n(n-3)/2 + P & \textrm{if $n>5, n\not=3k$}\\
n(n-3)/2 + P + $n/3[n+3/6]$ & \textrm{if $n>5, n=3k$}
\end{array} \right.
\end{displaymath}

P= the number of unordered pairs (a,b), where a,b are 2 parallel
diagonals in a regular n-gon.  It is not hard to find a precise
formula for P:  for every direction w, we have to count how many
diagonals are on this direction.  If there are z=z(w)
diagonals,then  ${z \choose 2}$ is the contribution of w.  Any
direction is given by the diagonals which pass through a fixed
vertex A.  We still have 3 more directions, given by AB, AC and
BC, where AB and AC are edges of the n-gon.

\vspace{10pt}

Let A be a maximal set of vertices of $G(n)$, such that the
distance between any 2 vertices is greater than 4. In this case
the small graphs of distance at most 2 from every vertex of A are
disjointed. They give independent non-zero cycles in $B(n)$ (see
the proof of Theorem 2).

So, the rank of the first homology is at least M times the
cardinal of A.

Theorem: The cardinal of A is at least (n-9)! ,if n is big enough.

\vspace{6pt}

Proof: A is maximal. so for any vertex x from G(n)-A , there is a
vertex b(x) from A such that the distance from x to b(x) is at
most 4. Otherwise, we could add x to A, in contradiction with the
maximality of A.

 \vspace{6pt}
Let's say that A is smaller than (n-9)!. Then, G(n)-A is greater
than n! - (n-9)! = (n-9)!*P(n), where P is a polynomial of degree
9. b is a function from a set with at least (n-9)!*P(n) elements,
to a set A of at most (n-9)! elements. So, there is an element y
from A, covered by at least P(n) elements from pre-image. So there
are at least P(n) elements at the distance at most 4 from y. The
degree of every vertex is n(n-3)/2, so at most $n^2$. There are at
most $n^4$ elements at the distance 2 from y. There are at most
$n^6$ elements at the distance 3 from y. There are at most $n^8$
elements at the distance 2 from y. Totally, there are at most
Q(n)= $n^2$ + $n^4$ + $n^6$ + $n^8$ vertices at distance at most 4
from x. For n big enough, P $>$ Q. Contradiction.

\vspace{10pt}
 \textbf{{\large The sketch of the proof for the lower bound M.}}

Definition. Let $G(n,k)$ and $B(n,k)$ be the graph and the
sub-complex of $B(n)$ generated by the vertices of distance at
most k from the vertex represented by the identity.

\vspace{10pt}

\textbf{{\large Theorem 2.}} The inclusion i:from $B(n,k)$ to
$B(n)$ gives a monomorphism in the 1-homology with real
coefficients, if k is smaller than n-3.

\vspace{5pt}

Note: the condition "k is smaller than n-3" is essential. G(n) is
a graph divided in levels but it natural metric. there are n-2
levels. All vertices from the first n-3 levels can be represented
by "permutations of type k" .A permutation of type k is a product
of k non-intersecting involutions ("diagonals") from T .

 Proof. Let a be a formal sum of edges  which satisfy
$i_{*}$(a)=0. So, a is a sum of boundaries= a sum of boundaries
from $B(n,k)$ + a sum of boundaries outside $B(n,k)$. The support
of a is in $B(n,k)$, so the last sum , let's call it A, is zero. A
is a sum of partial sums, given by a couple of non-intersecting
diagonals. The sum of the corresponding edges can be recovered by
boundaries inside $B(n,k)$, so a is zero in the homology of
$B(n,k)$.   So, the kernel of $i_{*}$ is zero.

\includegraphics [scale=2.0]{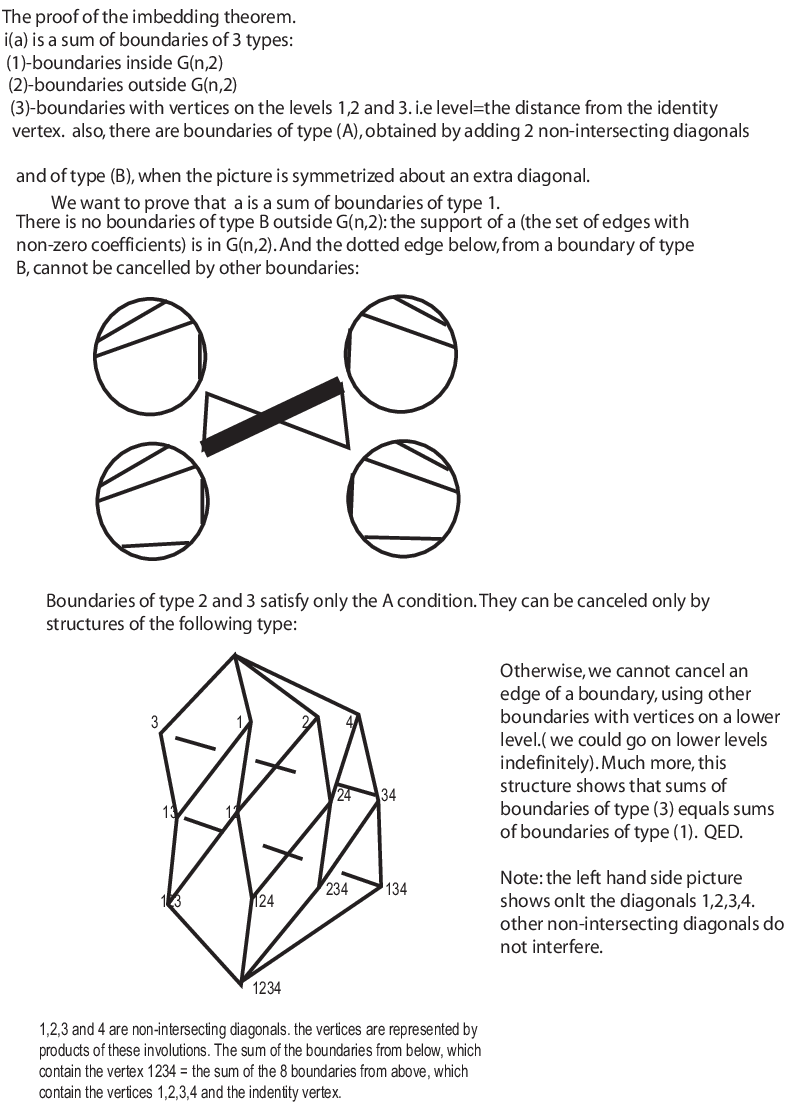}

\vspace{6pt}

\textbf{{\large Theorem 3.}} The rank of $H_{1}($B(n,2)$)$ is M.

We need "the theory of fluxes". from graph theory. For every
oriented graph G, let A be the real vector space of basis the
oriented edges . The vectors are formal sums of these edges, with
real coefficients. A flux is a formal sum which satisfy the
following property : for every vertex v, the sum of the
coefficients of the edges which enter in v is equal to the sum of
the coefficients of the edges which leave  v. The vector space of
all fluxes has the dimension m-n+1.[Bollobas pp.51-54], where m is
the number of edges and n is the number of vertices of G. This
formula gives us the opportunity to compute the rank as the
difference between 2 dimensions of vectors spaces (the kernel and
the imagine), the dimensions being separately computed.

\vspace{6pt}

The kernel is the space of fluxes. The image is studied using the
properties of the graph $G(n,2)$. Using a similar proof as above,
we can restrict to the complex $A(n,2)$, induced in $G(n,2)$ by
all vertices except the vertex represented by identity.

Notation: $(M-N+1)_{T}$ = the dimension of the space of fluxes in
the graph T, of M edges and N vertices.

The rank of $H_1(A(n,2))$ = $(M-N+1)_{A(n,2)}$ - dim(V U
W)=$(M-N+1)_{A(n,2)}$ - dim(V) - dim(W) + dim(V $\cap$ W), where

\par\noindent

V = the space generated by boundaries which are inside $A(n,2)$
(called "relations").

\par\noindent
W=the space of fluxes generated by boundaries outside $A(n,2)$,
which contain the identity vertex.
\par\noindent

M=the number of edges of $A(n,2)$= $n(n-3)/2 * [n(n-3)/2 -1]$

N=the number of vertices of $A(n,2)$=n(n-3)/2 + S, where S =1/3
${n-3 \choose 2}$ ${n+1 \choose 2}$ (the number of pairs of
non-intersecting diagonals).

dim(W)= ${m-n+1}_{W}$ = S - n(n-3)/2 +1

dim( V $\cap$ W)=the dimension of the space E generated by some
special cycles of length 3 (supported by lines which form an
equilateral triangle).

dim(V)=S - the number of pairs of parallel diagonals - A. A is the
number of independent relations between cycles of length 3,
supported by lines which form an equilateral triangle.  So it is
equal to the number of these cycles minus the dimension of the
space generated by them in W. Fortunately, we do not need to
compute this last number.  A + dim( V $\cap$ W) = the number of
equilateral triangles supported by 3 diagonals of an n-gon. We
apply the formulas above and we compute the announced result on M.

\includegraphics [scale=2.0]{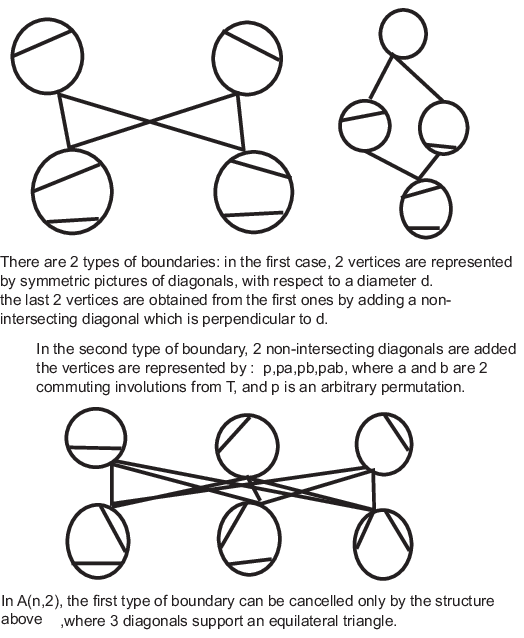}
\par\noindent
\includegraphics [scale=2.0]{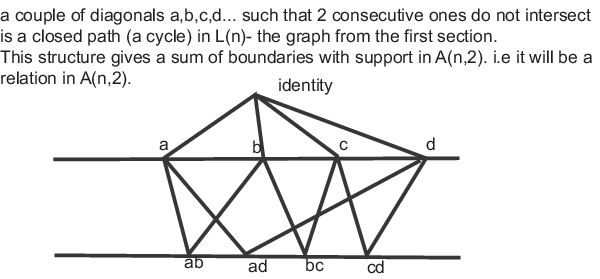}
\par\noindent
The result is a consequence of a combinatorial analysis on the
graph $G(n,2)$. It is fine enough to be carefully performed, but
it is elementary: it is a study of formal sums of n-gons with 1 or
2 diagonals, with real coefficients.

\par\noindent

\vspace{10pt}

\textbf{\texttt{Notes about the graphs G and $G(n,2)$:}
}diam(G)=n-2 (proof by induction).  Every vertex has degree
n(n-3)/2.  For any vertex v, take a permutation s from v, which is
a class modulo $D_n$. Then sd, where d runs in DIAG, represents
the neighbors of v.

2 vertices are joined by an edge in only 2 situations: "a
non-intersecting diagonal is added."  Or the picture is
symmetrized about a diameter perpendicular to a diagonal d, and d
disappears. More precisely: in these vertices, there are 2
permutations represented by a product of a couple of permutations
from DIAG, and there are only 2 situations when the vertices are
joined.

There are 2 types of good cycles of length 4:  when the opposite
vertices are on different levels, or when the opposite vertices
are on the same level. The first situation is possible when we add
2 non-intersecting diagonals.In the second situation, the 4
vertices are like:  $(a,b,ac,bc)$, where b is the symmetric of a
with respect to the diameter perpendicular on c (called the
symmetry diagonal). a and c are non-intersecting diagonals.

So the good cycles from V are linearly independent (because the
edges are met only one time, in these structures), except for one
situation: if there are a,b,c whose lines  form an equilateral
triangle, the 3 good cycles whose symmetry diagonals are a,b,c are
not independent:  their sum is a vector of W.  But we can have a
couple of these triads, and the sum can be zero (so we have a
relation).  This explains the computations of the dimensions of
vector spaces above.


\begin{thebibliography}{40}
\bibitem{morava} J.Morava.math.AT/0109086 talk Gdansk conference'01, pp.5
\bibitem{kont} M.Kontsevich: Feynman diagrams and low-dimensional topology. First European Congress of Math. Vol II(Paris, 1992), p.101
\bibitem{conf}Fadell, Edward R.; Husseini, Sufian Y. Geometry and topology of configuration spaces. Springer Monographs in Mathematics. Springer-Verlag, Berlin,
2001.
\bibitem{devadoss}S.Devadoss: Tessellations of moduli spaces and the mosaic operad. Homotopy invariant algebraic structures, 91-114 Contemp.Math. 239
\bibitem{barad}G.Barad: math.AG/0204003 The Homology of the real moduli spaces $\overline{M}_{0,n}$.
\bibitem{ceyhan}O.Ceyhan.Moduli of pointed real curves of genus 0.math.AG/0207058
\bibitem{davis}M.Davis, R.Scott, T.Januszkiewicz: Fundamental Groups of Blow-ups. Adv. in Math.to appear.
\bibitem{ka} M.Kapranov, Chow quotients of Grassmannians, I. Adv.in Soviet Math 16(1993), 29-111
\bibitem{ba} W.Ballmann, J.Swiatkowski: On $L^{2}$ -Cohomology and Property (T) for automorphism groups of polyhedral cell complexes. GAFA vol.7 (1997)
\bibitem{graph}Doob,Sachs,Cvetkovic: Spectra of Graphs.Theory and Appl. A.P.1980
\bibitem{bella}B.Bollobas: Modern Graph Theory. Springer'1998
\bibitem{s} \'Swi\c atkowski,J. Some infinite groups generated by involutions have Kazhdan's property (T).Forum Math(2001),6
\bibitem{chung}Chung, Fan R. K. Spectral graph theory. CBMS Regional Conference Series in Mathematics, 92
\bibitem{schreier} Gross, Jonathan L. Every connected regular graph of even degree is a Schreier coset graph. J. Combinatorial Theory Ser. B 22 (1977), no. 3
\bibitem{zi}G.Ziegler: Lectures on Polytopes. Springer'1995
\bibitem{zuk}\.Zuk, Andrzej La propriété (T) de Kazhdan pour les groupes agissant sur les polyèdres. (French) [Kazhdan's property (T) for groups acting on polyhedra] C. R. Acad. Sci. Paris Sér. I Math. 323 (1996), no. 5
\bibitem{loday}Loday, Jean-Louis; Ronco, María O. Une dualité entre simplexes standards et polytopes de Stasheff. (French) [A duality between standard simplices and Stasheff polytopes] C. R. Acad. Sci. Paris Sér. I Math. 333 (2001), no. 2
\bibitem{mu} Munkres: Elements of Algebraic Topology. Addison-Wesley'1984
\bibitem{phong}E.D'Hoker, D.H.Phong: Seiberg-Witten Theory and Calogero-Moser Systems. Progress of Theoretical Physics Suppl. No. 135, 1999. pp.75-93
\bibitem{vic}V.Reiner: Non-crossing partitions for classical reflection groups. Discrete Math. 177(1997) 195-222
\bibitem{ceyhan}Ceyhan,O.Moduli of pointed real curves of genus 0. math.AG/0207058
\bibitem{stanley}A polytope related to empirical distributions, plane trees, parking functions, and the associahedron. Discrete Comput. Geom. 27 (2002), no. 4
\bibitem{manin}Kontsevich, M.; Manin, Yu. Quantum cohomology of a product. With an appendix by R. Kaufmann. Invent. Math. 124 (1996), no. 1-3,
\end{thebibliography}
\end{document}